\documentclass{amsart}

\usepackage{parskip}

\usepackage{amsmath}
\usepackage{amsfonts}
\usepackage{amssymb}
\usepackage{latexsym}

\newtheorem{theorem}{Theorem}[section]

\def\N{{\Bbb N}}

\def\R{{\Bbb R}}
\def\C{{\Bbb C}}

\def\cC{{\mathcal C}}
\def\cD{{\mathcal D}}
\def\cH{{\mathcal H}}
\def\cL{{\mathcal L}}

\def\dim{\hbox{dim}}

\def\n*{\N_*^{m}}
\def\1{{\text{\bf 1}}}

\def\n{\vert\vert}
\def\eps{\epsilon}
\def\hol{{\mathcal O}}

\def\comp{\subset\subset}

\begin{document}

\title{ The H-principle and Pseudoconcave CR Manifolds }
\author{ C.~Denson Hill and Egmont Porten}
\date{\today}
\begin{abstract} 
The H-principle, which is the analogue, for CR manifolds, of the
classical Hartogs principle in several complex variables, is known
to be valid in the small on a pseudoconcave CR manifold of any
codimension. However it fails in the large, as has been shown by
the counterexample found in [HN1]. Hence there is an underlying
obstruction to the global H-principle on a pseudoconcave CR manifold.
The purpose of this note is to take the first steps toward a deeper
understanding of this obstruction.
\end{abstract}

\maketitle

\section{Introduction}
\hspace*{1cm}
Consider a smooth CR manifold $M$ of type $(n,k)$ which is a generic CR submanifold
of a complex manifold $X$. For a detailed discussion of all basic concepts concerning
CR manifolds, see any one of \cite{HN1}, \cite{HN2}, \cite{HN3}, \cite{HN4}.
Here $n$ is the CR dimension, $k$ is the CR codimension,
so $\dim_{\R} M=2n+k$, and $M$ being generic in $X$ implies that $\dim_\C X=n+k$.
When $k=0$, $M$ is just an $n$-dimensional complex manifold, and we take $X=M$.
In this note we shall always be considering the following situation: We assume
that {\it $\Omega$ is a connected open set in $M$, $K$ is a compact set in $\Omega$, 
and $\Omega\backslash K$ is connected and nonvoid}. $CR(D)$ denotes
the algebra of continuous CR functions in the open set $D\subset M$, i.e., functions
which satisfy the tangential Cauchy-Riemann equations on $M$ in the sense of
distributions. When $k=0$ the CR  functions are simply the holomorphic functions,
and we write $\hol(D)$ instead of $CR(D)$.

In this situation we say that {\it the H-principle holds for the pair} $(\Omega,K)$
iff the restriction map
\begin{equation}\label{1.1}
CR(\Omega)\longrightarrow CR(\Omega\backslash K)
\end{equation}
is surjective. Here $K$ should be thought of as a ``hole'' which can always be ``filled''
by CR functions. If the CR functions on $M$ enjoy the weak unique continuation property
(i.e., vanishing in an open set implies vanishing everywhere), then 
the restriction map (\ref{1.1}) is injective, so in that case the validity of the H-principle 
actually means that there is an isomorphism $CR(\Omega)\cong CR(\Omega\backslash K)$ of 
algebras. 

When $k=0$, (\ref{1.1}) becomes $\hol(\Omega)\longrightarrow\hol(\Omega\backslash K)$.
If this map is surjective, then it is an isomorphism, and the hole $K$ in $\Omega$ can
always be filled by holomorphic extension; this is 
the {\it classical Hartogs-principle (H-principle)} in several complex variables. The
classical Hartogs-principle is valid for any pair $(\Omega,K)$, provided the 
$n$-dimensional complex manifold $M$ is $(n-2)$-complete (in particular when $n>1$,
and $M$ is $\C^n$, or $M$ is Stein), for example see \cite{AH}.

When $k\geq 1$ the H-principle clearly fails, as is seen by easy counterexamples,
{\it unless} $M$ has some amount of pseudoconcavity (perhaps in a weak sense).
A CR manifold $M$
is said to be $q$-pseudoconcave provided that the Levi form $\cL_{x_0}(\xi,\cdot)$ of $M$
at $x_0$ has at least $q$ negative eigenvalues in each nonzero characteristic conormal
direction $\xi$ at $x_0$, for every $x_0\in M$. 

A $1$-pseudoconcave CR manifold is simply called pseudoconcave. A pseudoconcave CR manifold 
$M$ enjoys the property that to any connected
open set $D\subset M$, one can associate a connected open set $\tilde{D}\subset X$, with
$D=\tilde{D}\cap M$, such that the restriction map $\hol(\tilde{D})\rightarrow CR(D)$ is
an isomorphism of algebras (see \cite{BP}, \cite{NV}, \cite{HN5}); 
so in particular, the CR functions on $M$ have the strong unique
continuation property (i.e., all derivatives vanishing at a point implies vanishing
everywhere), and are as smooth as $M$ is. So in this situation the CR functions 
behave very much like holomorphic functions. Thus on a pseudoconcave CR manifold $M$,
it seems to be a natural question to ask: {\it When does the H-principle hold for certain pairs
$(\Omega,K)$?}

Indeed it was shown in \cite{He} that the H-principle is valid on a smooth pseudoconcave
CR manifold $M$ for any pair $(\Omega,K)$, {\it provided that $\Omega$ is sufficiently
small}. For real analytic $M$ and real analytic CR functions, this result has been extended
\cite{HeM} to the weakly 1-pseudoconcave case. See also \cite{Na} for an earlier result,
in a special case.

On the other hand a very simple counterexample was found in \cite{HN1} for the situation
where $K$ becomes too big: Take $M$ to be the pseudoconcave CR hypersurface of type (2,1) 
in $\C^3$ given by
\begin{equation}\label{1.2}
M:|z_1|^2+|z_2|^2=1+|z_3|^2,
\end{equation}
and $K=S^3$ to be $M\cap\{z_3=0\}$. By considering the CR function $1/z_3$ on $M\backslash S^3$,
one sees that the H-principle {\it cannot hold} for any pair $(\Omega,S^3)$. Thus there is an
underlying obstruction to the global H-principle for pseudoconcave $M$. The present note is
intended as a first step toward a deeper understanding of this obstruction.

We shall show below that example (\ref{1.2}) has some interesting features:
\begin{description}
\item[(A)] 
If $S^3$ is replaced by $K_\omega=S^3\backslash\omega$, where $\omega$ is an arbitrarily
small non-empty open set on $S^3$, then the H-principle {\it holds} for any pair 
$(\Omega,K_\omega)$.
\item[(B)] 
If $S^3$ is replaced by $K=\tilde{S}^3$, where $\tilde{S}^3\subset M$ is a small, but randomly
chosen smooth perturbation of $S^3$, then the H-principle {\it holds} for any pair 
$(\Omega,\tilde{S}^3)$.
\end{description}
We also obtain that both $CR(M\backslash K_\omega)$ and $CR(M\backslash\tilde{S}^3)$ are isomorphic 
to $\hol(\C^3)$, the space of entire functions on $\C^3$.

Note that in both (A) and (B) above, the hole which is filled by the H-principle is ``thin'', in the
sense that both $K_\omega$ and $\tilde{S}^3$ have codimension $2$ in $M$. One can ask about trying
to fill ``thick'' holes, by which we mean that $K$ should be a compact set having non-empty interior.
In this respect (A) and (B) behave quite differently. 
\begin{description}
\item[(C)] 
Given a pair $(\Omega,K_\omega)$, there exists an $\eps=\eps(\Omega,K_\omega)>0$ such that
the H-principle {\it holds} for any $(\Omega,K_\omega^\eps)$, where  $K_\omega^\eps$ denotes the closed
$\eps$-neighborhood of $K_\omega$ in $M$.
\item[(D)]
But if $\tilde{S}^3$ in (B) is replaced by any slender compact neighborhood $\tilde{K}$ 
of $\tilde{S}^3$ in $M$ which contains the original $S^3$,
then the H-principle {\it cannot hold} for any pair $(\Omega,\tilde{K})$, as is obvious
because of the original example.
\end{description}
We should also draw the readers attention to the interesting papers \cite{L}, \cite{LL1}, 
\cite{LL2}, \cite{B},
in which global results are obtained under additional assumptions, mainly on the degree of
pseudoconcavity.

The first author would like to acknowledge the support of the Mathematics Institute of Humboldt
University in Berlin, and in particular the kind hospitality of Professor J\"urgen Leiterer.

\section{Elementary use of analytic discs.}
\setcounter{equation}{0}
\hspace*{1cm}
In order to motivate the general results in section 3, we first discuss points (A) and (C)
for example (\ref{1.2}).

We consider parameters $(r,s)$ with $0\leq r<\infty$ and $s\in S^3$, where $S^3$  denotes the
unit sphere in $\C^2$ with coordinates $(z_1,z_2)$. Let $\Delta_{r,s}$ denote the subset of
$\C^3$ defined by
\begin{equation}\label{2.1}
\Delta_{r,s}=\{([1+r^2]^{1/2}s,z_3)\big| |z_3|<r\},
\end{equation}
and $\partial\Delta_{r,s}=\overline{\Delta}_{r,s}\backslash\Delta_{r,s}$ denote its boundary.
Note that the $\Delta_{r,s}$ give a foliation of the exterior $\{|z_1|^2+|z_2|^2>1+|z_3|^2\}$
of our $M$ given in (\ref{1.2}), the boundaries $\partial\Delta_{r,s}$ lie on $M$, and shrink 
to the points $(s,0)$ when $r=0$.

In order to verify statement (A) we choose an open set $V$ on $S^3$ such that $K_\omega\subset
V\comp S^3\cap\Omega$, and $r_2>0$ sufficiently small so that
\begin{equation}\label{2.2}
\bigcup_{0<r<r_2,s\in V}\partial\Delta_{r,s}\subset\Omega.
\end{equation}
Fix $(\hat{r},\hat{s})\in(0,r_2)\times V$ and choose in $V$ a continuous path $\gamma(t),0\leq
t\leq 1$, connecting $\hat{s}$ with some $\tilde{s}\in V\backslash K_\omega$. Then for $0\leq
t\leq 1$ the boundaries  $\partial\Delta_{\hat{r},\gamma(t)}\subset\Omega\backslash K_\omega$.
But by shrinking $\hat{r}$ to zero, the boundary $\partial\Delta_{\hat{r},s}$ can be
contracted to the point $\tilde{s}\in\Omega\backslash K_\omega$. As was mentioned in the
introduction, the CR functions on $\Omega\backslash K_\omega$ have holomorphic extensions
to an ambient neighborhood of $\Omega\backslash K_\omega$ in $\C^3$. By the Kontinuit\"atssatz
we therefore obtain holomorphic extension to
\begin{equation}\label{2.3}
\bigcup_{0<r<r_2,s\in V}\Delta_{r,s}\subset\C^3,
\end{equation}
which is an open set attached to an open neighborhood of $K_\omega$ in $M$, from the exterior.
By the pseudoconcavity of $M$ we obtain holomorphic extension across $K_\omega$, proving (A).

In order to prove (C) for the pair $(\Omega,K^\eps_\omega)$ it suffices to modify the above
argument as follows: With $V$ and $r_2$ chosen exactly as above, we first remark that the 
same proof works if $K_\omega$ is replaced by $K_{\omega,\delta}$, where $K_{\omega,\delta}$
is the closed $\delta$-neighborhood of $K_\omega$ in $S^3$, and $\delta>0$ is sufficiently
small. In fact the envelope of holomorphy of the set (\ref{2.3}) still contains 
$K_{\omega,\delta}$. Next we observe that if $r_1>0$ is chosen sufficiently small, then the
envelope of holomorphy of
\begin{equation}\label{2.4}
\bigcup_{r_1<r<r_2,s\in V}\partial\Delta_{r,s}\subset\Omega,
\end{equation}
contains a neighborhood $W$ of $K_{\omega,\delta}$ in $M$. Hence it suffices to choose $\eps>0$
so small that every function in $\hol(\Omega\backslash K^\eps_\omega)$ extends holomorphically
to the set (\ref{2.4}) and that $K^\eps_\omega\subset W$. This completes the proof of point (C).

A direct proof of (B) would be considerably more involved. However we obtain it as a special
case of Theorem 3.1 in the next section.

\section{Removing thin singularities}
\setcounter{equation}{0}
\hspace*{1cm}
In this section we return to the general situation of a CR manifold $M$ of type $(n,k)$, 
as in the beginning of the Introduction, and we study the H-principle for thin holes $K$.

Consider a point $x_0\in M$ and a local system $z=(z_1,z_2,\ldots,z_{n+k})$ of 
holomorphic coordinates for $X$, centered at $x_0$. Then $M$ is locally defined by
\begin{equation} \label{3.1}
   M = \{ \rho_1(z)=0, \rho_2(z)=0,\ldots, \rho_k(z)=0\}
\end{equation}
where the $\rho_i$ are smooth real valued functions, defined in a neighborhood of $x_0$. 
The fact that $M$ is generic in $X$ means that 
               $\overline{\partial}\rho_1,\overline{\partial}\rho_2,\ldots,\overline{\partial}\rho_k$ 
are linearly independent at $x_0$. 
The holomorphic tangent space to $M$ at $x_0$ is defined by
\begin{equation} \label{3.2}
   H_{x_0} M = \{ w \in \C^{n+k} \mid \sum_{j=1}^{n+k} \frac{\partial \rho_i(x_0)}{\partial z_j} w_j =0; 
                                                                                             i=1,2,\ldots,k \}.
\end{equation}
The Levi form ${\mathcal L}_{x_0} (\xi,\cdot)$ of $M$ at $x_0$, in the characteristic covector direction 
$\xi=(\xi_1,\xi_2,\ldots,\xi_k)$, can be written as
\begin{equation} \label{3.3}
   \sum_{i=1}^k \sum_{j,l=1}^{n+k} \xi_i 
              \frac{\partial^2 \rho_i(x_0)}{\partial z_j \partial \overline{z}_l} w_j \overline{w}_l .
\end{equation}
It is a hermitian form on $H_{x_0}M \cong \C^n$, and it is useful to write 
$T_{x_0}M \cong \C^n\times\R^k$.

Next let $N$ be a smooth closed submanifold of $M$ having real codimension 2 in $M$. 
Then $N$ is locally defined by
\begin{equation} \label{3.4}
   N = \{ \rho_1(z)=0,\ldots,\rho_k(z)=0,\rho_{k+1}(z)=0,\rho_{k+2}(z)=0 \}
\end{equation}
in a neighborhood of $x_0 \in N$, with $d\rho_1,\ldots,d\rho_k,d\rho_{k+1},d\rho_{k+2}$ linearly
independent at $x_0$. Since we are not assuming that $N$ is a CR manifold, the
\begin{equation} \label{3.5}
   \mbox{span} \{ \overline{\partial}\rho_1,\ldots,\overline{\partial}\rho_k,
   \overline{\partial}\rho_{k+1},
    \overline{\partial}\rho_{k+2} \}
\end{equation}
at $x_0$ may have dimension $k+2$, $k+1$ or $k$. Thus the 3 possibilities are $T_{x_0}N \cong \C^{n-2}\times\R^{k+2}$,
$T_{x_0}N \cong \C^{n-1}\times\R^{k}$ or $T_{x_0}N \cong \C^{n}\times\R^{k-2}$. In the case where 
$T_{x_0}N \cong \C^{n-2}\times\R^{k+2}$, $N$ is said to be {\it generic} at $x_0$.

In the theorem below we consider a smooth CR manifold $M$ of type $(n,k)$, as in the Introduction, 
and a smooth closed connected submanifold $N$ of $M$.

\begin{theorem}\label{cd2}
Assume $M$ is pseudoconcave, $N$ has real codimension 2 in $M$, and let $K$ be a compact subset of $N$.
Then the H-principle is valid for any pair $(\Omega,K)$  provided one of the following conditions 
holds:

\noindent
(a) $K$ is a proper subset of $N$, or \\
(b) $K=N$ and there is a point $x_0\in N$ at which $N$ is generic.
\end{theorem}

Let us interpret Theorem \ref{cd2} in terms of the pseudoconcave $M$ 
of type $(2,1)$ in example (\ref{1.2}):
If we choose $N=S^3$, then from (a) we recover our statement (A). 
However (a) gives much more; namely it
yields an analogous result for any smooth 3-dimensional submanifold $N$ of $M$.

If we choose $N=\tilde{S}^3$, then from (b) we obtain our remaining statement (B). 
This is because a randomly chosen deformation of $S^3$ will contain many generic 
points. However (b) gives much more; namely any codimension 2
counterexample $K=N$ has to be a CR manifold of type $(1,1)$.

{\bf Proof of Theorem \ref{cd2}.} Because of the uniqueness of CR extension on $M$, 
it will suffice to show 
that we have CR extension across any given point $x_0\in K$. 
As was mentioned in the Introduction, any given 
$f\in \mbox{\it CR} (\Omega\backslash K)$ has a holomorphic extension, again denoted by $f$, to an open set 
$\tilde{D}$ in $X$ with $\Omega\backslash K = \tilde{D}\cap M$. 
By \cite[Theorem 4]{MP1}, see also \cite{CS}, \cite{J}, and \cite{M}, there is an open truncated wedge
$W \subset X$ attached to an open neighborhood $U$ of $x_0$ in $M$, an open set $\hat{D}$ of $X$ with 
$\hat{D}\subset\tilde{D}$ and with $\Omega\backslash K = \hat{D}\cap M$, and an $f_w \in \hol(W)$ such that 
$f_w \mid_{W\cap\hat{D}} = f \mid_{W\cap\hat{D}}$. Here $W$ and $\hat{D}$ do not depend on $f$. 

We have to justify that the theorem being used here applies under either hypothesis 
(a) or (b). First we observe that in both cases $M$ and $M\backslash K$ are globally
minimal. Indeed, being pseudoconcave, $M$ is even minimal in each of its points. This can
be seen as follows: In each $p\in M$, the brackets of the vectorfields tangent to $H$
span the whole tangent space $T_p M$. Since these brackets are tangent to the local
CR orbit of $M$ in $p$, the local CR orbit has to be open 
(see \cite{S}, \cite{J}, \cite{MP1} for detailed information on CR orbits).
In case (b) we also need $n\geq 2$, which is an obvious consequence of
pseudoconcavity.

In terms of the local holomorphic 
coordinates $z$ in $\C^{n+k}$, centered at $x_0$, $W$ can be chosen as follows: For some open truncated cone
$C\subset\C^{n+k}$, with vertex at the origin, $W = U+C = \{ z+c \mid z\in U, c\in C\}$. However associated 
to the open set $U$ in $M$ there is an open set $\tilde{U}$ in $\C^{n+k}$, with $U = \tilde{U} \cap M$, such 
that $\hol(\tilde{U}) \cong \mbox{\it CR} (U)$. 
For any vector $c\in C$ we define the rigid motion translates of $U$
by $U_c = U+c$, and note that $U_c\subset W$. Note that $\hol(\tilde{U}_c) \cong \mbox{\it CR}(U_c)$ where
$\tilde{U}_c = \tilde{U}+c$. By choosing $\mid c\mid$ sufficiently small, we have $x_0 \in \tilde{U}_c$.
This gives the desired extension of our original $f$ to a neighborhood of $x_0$, and completes the proof of
the theorem.

We end this section with a second theorem which shows that for pseudoconcave CR manifolds $M$ the H-principle
can only fail for holes $K$ having codimension at most 2. Let $\cH^s(K)$ denote s-dimensional Hausdorff measure 
of $K$.

\begin{theorem}\label{thin}
Assume $M$ as in Theorem \ref{cd2}, and let $K\subset M$ be any compact subset with $\cH^{2n+k-2}(K)=0$.
Then the H-principle is valid for any pair $(\Omega,K)$.
\end{theorem}

{\bf Proof.} The proof is the same as that of Theorem \ref{cd2}, except that we use 
\cite[Theorem 1.1]{MP2}, see also \cite{LS}, \cite{CS}, to obtain
holomorphic extension into a wedge $W$ attached to a neighborhood of $x_0\in K$.

\section{Revisiting the example}
\setcounter{equation}{0}
\hspace*{1cm}
We now return to the discussion of the specific $M$ given by (\ref{1.2}).
We know from (B) that for a pair $(M,\tilde{S}^3)$ the H-principle is {\it valid} for
most small deformations $\tilde{S}^3$ of $S^3$. In the next theorem we characterize
precisely those small deformations for which the H-principle for $(M,\tilde{S}^3)$
{\it fails}.

\begin{theorem}\label{nesu}
Let $\tilde{S}^3$ be a sufficiently small $\cC^2$-deformation of $S^3$. Then the 
H-principle for $(M,\tilde{S}^3)$ fails if and only if $\tilde{S}^3$ is the 
intersection of $M$ with a complex hypersurface $Y$ in $\C^3$.
\end{theorem}

{\bf Proof.}
If there is such a $Y$, then there is an entire function $g$ on $\C^3$ such
that $Y=\{z\in\C^3|g(z)=0\}$. Then $1/g$ is a CR function on $M\backslash\tilde{S}^3$
which has no CR extension across $\tilde{S}^3$, so the H-principle fails.

On the other hand, if the H-principle fails, then we know that $\tilde{S}^3$ must be
a CR manifold of type (1,1). Let $\pi:\C^3\rightarrow\C^2$ and $\pi_3:\C^3\rightarrow\C$
denote the holomorphic projections onto the $(z_1,z_2)$-plane and the $z_3$-axis,
respectively. Since the deformation is $\cC^2$-small, we obtain $\tilde{S}^3$ as 
a graph of the CR function $\phi=\pi_3\circ(\pi|_{\tilde{S}^3})^{-1}$ over the
stricly convex 3-dimensional hypersurface $\pi(\tilde{S}^3)$, which bounds a domain
$G$ in $\C^2$. It is well known that $\phi$ has an extension $\tilde{\phi}\in
\hol(G)\cap\cC^2(\overline{G})$. The graph of $\tilde{\phi}$ over $G$ defines
a smooth complex hypersurface $\tilde{Y}$ bounded by $\tilde{S}^3$. Note that
$\tilde{Y}$ is contained in the interior of the domain 
$\cD=\{|z_1|^2+|z_2|^2<1+|z_3|^2\}$ and is transversal to its boundary $M$.

Next we show that any $f\in\mbox{\it CR}(M\backslash\tilde{S}^3)$ has a 
holomorphic extension $\tilde{f}$ to $\cD\backslash\tilde{Y}$. Let $\eps>0$
be given. By Runge approximation there is a polynomial $P_\eps\in\C[z_1,z_2]$
such that $\max_{\overline{G}}|P-\tilde{\phi}|<\eps$. Then the complex hypersurface
\begin{equation}\label{4.1}
Y_{\eps,c}=\{z_3=P_\eps(z_1,z_2)+c\}
\end{equation}
does not intersect $\tilde{Y}$ provided $|c|>\eps$. Since each $Y_{\eps,c}$ is 
Stein, the H-principle implies that $f$ has a holomorphic extension to
$\overline{\cD}\cap Y_{\eps,c}$ for $|c|>\eps$, and we obtain the holomorphic
extension of $f$ to
\begin{equation}\label{4.2}
\cD_\eps={\mathcal D}\cap\{|z_3-P_\eps(z_1,z_2)|>\eps\}.
\end{equation}
Taking $\eps\searrow 0$ we see that $f$ has a holomorphic extension to 
$\cD\backslash\tilde{Y}$.

By \cite{D}, see also \cite{P}, there are only two possibilities for 
the envelope of holomorphy $\Sigma$ of ${\mathcal D}\backslash \tilde{Y}$:
Either (1) $\Sigma = \C^3$ or else (2) there is a complex hypersurface $Y$ of $\C^3$ 
satisfying $\tilde{Y} = Y\cap\mathcal D$ such that $\Sigma=\C^3\backslash Y$. 
Here we have used that the envelope of holomorphy of $\mathcal D$
is all of $\C^3$. Possibility (1) is ruled out by our hypothesis. Thus it remains only to 
verify that $Y\cap M = \tilde{S}^3$. Obviously $\tilde{S}^3 \subset Y\cap M$. 
But there are no points in $(M\backslash \tilde{S}^3)\cap Y$
because $M$ is pseudoconcave. This completes the proof of the theorem.

{\bf Remarks.} 
It follows from the argument above, that if the H-principle fails for $(M,\tilde{S}^3)$, 
where $\tilde{S}^3\subset M$ is a small $\cC^2$-deformation of $S^3$, then $\tilde{S}^3$ must 
be a {\it real analytic} CR manifold of type $(1,1)$.    

Let $\Omega$ be a given open set on $M$ containing $S^3$. Again consider $\tilde{S}^3
\subset\Omega$ to be a sufficiently small
$\cC^2$-deformation of $S^3$. If the H-principle fails for $(\Omega,\tilde{S}^3)$ 
then, once again $\tilde{S}^3$
must be a {\it real analytic} CR manifold of type $(1,1)$. 
This follows by a straightforward modification of the 
proof of Theorem \ref{nesu}, in which the role of $\C^3$ is replaced by the 
envelope of holomorphy of an appropriate 
domain, analogous to $D\backslash \tilde{Y}$.

\quad

\noindent%
C.~Denson Hill \\             
Department of Mathematics \\
Stony Brook University \\
Stony Brook, N.Y. 11794, USA \\
dhill@math.sunysb.edu \\  

\noindent%
Egmont Porten \\
Mathematisches Institut \\
Humboldt-Universit\"at zu Berlin \\
Rudower Chaussee 25 \\      
12489 Berlin, Germany \\
egmont@mathematik.hu-berlin.de          
             
\end{document}